# On the Developable Mannheim Offsets of Timelike Ruled Surfaces

**Mehmet Önder[*], H. Hüseyin Uğurlu[**]**

[*]Celal Bayar University, Faculty of Science and Arts, Department of Mathematics, Muradiye Campus, Muradiye, Manisa, Turkey. E-mail: mehmet.onder@cbu.edu.tr

[**]Gazi University, Faculty of Education, Department of Secondary Education Science and Mathematics Teaching, Mathematics Teaching Programme, Ankara, Turkey.
E-mail: hugurlu@gazi.edu.tr

**Abstract**
In this paper, using the classifications of timelike and spacelike ruled surfaces, we study the Mannheim offsets of timelike ruled surfaces in Minkowski 3-space. First, we define the Mannheim offsets of a timelike ruled surface by considering the Lorentzian casual character of the offset surface. We obtain that the Mannheim offsets of a timelike ruled surface may be timelike or spacelike. Furthermore, we characterize the developable of Mannheim offset of a timelike ruled surface by the derivative of the conical curvature $\kappa$ of the directing cone.

**MSC:** 53A10, 53B25, 53C50.

**Key words:** Mannheim offsets; Minkowski 3-space; timelike ruled surface.

## 1. Introduction

Ruled surfaces are the surfaces which are generated by moving a straight line continuously in the space and are one of the most important topics of differential geometry. These surfaces have an important role and many applications in the study of design problems in spatial mechanism, physics and Computer Aided Geometric Design (CAGD). Because of these positions of the ruled surfaces, many geometers have studied on them in Euclidean space and they have investigated the many properties of the ruled surfaces [5,11,12,13]. Furthermore, the differential geometry of the ruled surfaces in Minkowski space has been studied in the references [1,4,6,7]. Moreover, using the classification of the ruled surfaces, Uğurlu and Önder have given the Frenet frames, invariants and instantaneous rotation vectors of the timelike and spacelike ruled surfaces in Minkowski 3-space $IR_1^3$ [15,16].

Furthermore, in the plane, while a curve $\vec{\alpha}$ rolls on a straight line, the center of curvature of its point of contact describes a curve $\vec{\beta}$ which is the Mannheim partner of $\vec{\alpha}$. Mannheim partner curves have been studied by Liu and Wang in three dimensional Euclidean 3-space and Minkowski 3-space [8,18]. They have given the definition of Mannheim partner curves as follows: Let $C$ and $C^*$ be two space curves. $C$ is said to be a Mannheim partner curve of $C^*$ if there exists a one to one correspondence between their points such that the binormal vector of $C$ is the principal normal vector of $C^*$. They have shown that $C$ is Mannheim partner curve of $C^*$ if and only if
$$\frac{d\tau}{ds} = \frac{\kappa}{\lambda}(1+\lambda^2\tau^2),$$
where $\kappa$ and $\tau$ are the curvature and the torsion of the curve $C$, respectively, and $\lambda$ is a nonzero constant.

Considering the notion of Bertrand curves, Ravani and Ku have defined and studied the Bertrand offsets of the ruled surfaces [13]. By a similar way, the Mannheim offsets of ruled surfaces have been defined and studied by Orbay and et al in 3-dimensional Euclidean space $E^3$ [10].

In this paper, by considering the classifications of the ruled surfaces in Minkowski 3-space, we give the Mannheim offsets of the timelike ruled surfaces in Minkowski 3-space $IR_1^3$.

## 2. Preliminaries

The Minkowski 3-space $IR_1^3$ is the real vector space $IR^3$ provided with the standard flat metric given by
$$\langle , \rangle = -dx_1^2 + dx_2^2 + dx_3^2$$
where $(x_1, x_2, x_3)$ is a standard rectangular coordinate system of $IR_1^3$. An arbitrary vector $\vec{v} = (v_1, v_2, v_3)$ in $IR_1^3$ can have one of three Lorentzian causal characters; it can be spacelike if $\langle \vec{v}, \vec{v} \rangle > 0$ or $\vec{v} = 0$, timelike if $\langle \vec{v}, \vec{v} \rangle < 0$ and null (lightlike) if $\langle \vec{v}, \vec{v} \rangle = 0$ and $\vec{v} \neq 0$. Similarly, an arbitrary curve $\vec{\alpha} = \vec{\alpha}(s)$ can locally be spacelike, timelike or null (lightlike), if all of its velocity vectors $\vec{\alpha}'(s)$ are spacelike, timelike or null (lightlike), respectively. We say that a timelike vector is future pointing or past pointing if the first compound of the vector is positive or negative, respectively. The norm of the vector $\vec{v} = (v_1, v_2, v_3) \in IR_1^3$ is given by
$$\|\vec{v}\| = \sqrt{|\langle \vec{v}, \vec{v} \rangle|}.$$

For any vectors $\vec{x} = (x_1, x_2, x_3)$ and $\vec{y} = (y_1, y_2, y_3)$ in $IR_1^3$, Lorentzian vector product of $\vec{x}$ and $\vec{y}$ is defined by
$$\vec{x} \times \vec{y} = \begin{vmatrix} e_1 & -e_2 & -e_3 \\ x_1 & x_2 & x_3 \\ y_1 & y_2 & y_3 \end{vmatrix} = (x_2 y_3 - x_3 y_2, x_1 y_3 - x_3 y_1, x_2 y_1 - x_1 y_2).$$
(For details see [9]). The Lorentzian sphere and hyperbolic sphere of radius $r$ and center 0 in $IR_1^3$ are given by
$$S_1^2 = \{\vec{x} = (x_1, x_2, x_3) \in IR_1^3 : \langle \vec{x}, \vec{x} \rangle = r^2\}$$
and
$$H_0^2 = \{\vec{x} = (x_1, x_2, x_3) \in IR_1^3 : \langle \vec{x}, \vec{x} \rangle = -r^2\},$$
respectively [17].

**Definition 2.1.** *i) Hyperbolic angle:* Let $\vec{x}$ and $\vec{y}$ be future pointing (or past pointing) timelike vectors in $IR_1^3$. Then there is a unique real number $\theta \geq 0$ such that $<\vec{x}, \vec{y}> = -\|\vec{x}\|\|\vec{y}\| \cosh \theta$. This number is called the *hyperbolic angle* between the vectors $\vec{x}$ and $\vec{y}$.

*ii) Lorentzian timelike angle:* Let $\vec{x}$ be a spacelike vector and $\vec{y}$ be a timelike vector in $IR_1^3$. Then there is a unique real number $\theta \geq 0$ such that $<\vec{x}, \vec{y}> = \|\vec{x}\|\|\vec{y}\| \sinh \theta$. This number is called the *Lorentzian timelike angle* between the vectors $\vec{x}$ and $\vec{y}$.

**Definition 2.2.** A surface in the Minkowski 3-space $IR_1^3$ is called a timelike surface if the induced metric on the surface is a Lorentz metric and is called a spacelike surface if the induced metric on the surface is a positive definite Riemannian metric, i.e., the normal vector on the spacelike (timelike) surface is a timelike (spacelike) vector [2].

## 3. Differential Geometry of the Ruled Surfaces in Minkowski 3-space

Let $I$ be an open interval in the real line $IR$, $\vec{k} = \vec{k}(s)$ be a curve in $IR_1^3$ defined on $I$ and $\vec{q} = \vec{q}(s)$ be a unit direction vector of an oriented line in $IR_1^3$. Then we have the following parametrization for a ruled surface $N$

$$\varphi(s,v) = \vec{k}(s) + v\vec{q}(s). \tag{1}$$

The parametric $s$-curve of this surface is a straight line of the surface which is called ruling. For $v = 0$, the parametric $v$-curve of this surface is $\vec{k} = \vec{k}(s)$ which is called base curve or generating curve of the surface. In particular, if the direction of $\vec{q}$ is constant, the ruled surface is said to be cylindrical, and non-cylindrical otherwise.

The striction point on a ruled surface $N$ is the foot of the common normal between two consecutive rulings. The set of the striction points constitutes a curve $\vec{c} = \vec{c}(s)$ lying on the ruled surface and is called striction curve. The parametrization of the striction curve $\vec{c} = \vec{c}(s)$ on a ruled surface is given by

$$\vec{c}(s) = \vec{k}(s) - \frac{\langle d\vec{q}, d\vec{k} \rangle}{\langle d\vec{q}, d\vec{q} \rangle} \vec{q}. \tag{2}$$

So that, the base curve of the ruled surface is its striction curve if and only if $\langle d\vec{q}, d\vec{k} \rangle = 0$. Furthermore, the generator $\vec{q}$ of a developable ruled surface is tangent of its striction curve.

The distribution parameter (or drall) of the ruled surface in (1) is given as

$$\delta_q = \frac{\left| d\vec{k}, \vec{q}, d\vec{q} \right|}{\langle d\vec{q}, d\vec{q} \rangle} \tag{3}$$

(see [1,14]). If $\left| d\vec{k}, \vec{q}, d\vec{q} \right| = 0$, then the normal vectors are collinear at all points of the same ruling and at the nonsingular points of the surface $N$, the tangent planes are identical. Then we say that the tangent plane contacts the surface along a ruling. Such a ruling is called a *torsal* ruling. If $\left| d\vec{k}, \vec{q}, d\vec{q} \right| \neq 0$, then the tangent planes of the surface $N$ are distinct at all points of the same ruling which is called nontorsal [15,16].

**Definition 3.1.** A ruled surface whose all rulings are torsal is called a *developable ruled surface* in $IR_1^3$. The remaining ruled surfaces are called skew ruled surfaces [15,16].

**Theorem 3.1.** *In* $IR_1^3$, *a ruled surface is developable if and only if at all its points the distribution parameter* $\delta_q = 0$ [1,14,15,16].

For the unit normal vector $\vec{m}$ of the ruled surface $N$ we have $\vec{m} = \frac{\vec{\varphi}_s \times \vec{\varphi}_v}{\|\vec{\varphi}_s \times \vec{\varphi}_v\|}$.

So, at the points of a nontorsal ruling $s = s_1$ we have

$$\vec{a} = \lim_{v \to \infty} \vec{m}(s_1, v) = \frac{(d\vec{q}/ds) \times \vec{q}}{\|d\vec{q}/ds\|}.$$

The plane of the ruled surface $N$ which passes through its ruling $s_1$ and is perpendicular to the vector $\vec{a}$ is called the *asymptotic plane* $\alpha$. The tangent plane $\gamma$ passing through the ruling $s_1$ which is perpendicular to the asymptotic plane $\alpha$ is called the *central plane*. Its point of contact $C$ is *central point* of the ruling $s_1$. The straight lines which pass through point $C$ and are perpendicular to the planes $\alpha$ and $\gamma$ are called the *central tangent* and *central normal*, respectively.

Since the vectors $\vec{q}$, $d\vec{q}/ds$ and $\vec{a}$ are perpendicular, representation of the unit vector $\vec{h}$ of the central normal is given by

$$\vec{h} = \frac{d\vec{q}/ds}{\|d\vec{q}/ds\|}.$$

The orthonormal system $\{C; \vec{q}, \vec{h}, \vec{a}\}$ is called Frenet frame of the ruled surfaces $N$ such that $\vec{h} = \frac{d\vec{q}/ds}{\|d\vec{q}/ds\|}$ and $\vec{a} = \frac{(d\vec{q}/ds) \times \vec{q}}{\|d\vec{q}/ds\|}$ are the central normal and the asymptotic normal direction of $N$, respectively, and $C$ is the striction point.

Let now consider the ruled surface $N$ with non-null Frenet vectors and their non-null derivatives. According to the Lorentzian characters of the ruling and central normal, we can give the following classifications of the timelike or spacelike ruled surface $N$ as follows;

**i)** If the central normal vector $\vec{h}$ is spacelike and $\vec{q}$ is timelike, then the ruled surface $N$ is said to be of type $N_-^1$.

**ii)** If the central normal vector $\vec{h}$ and the ruling $\vec{q}$ are both spacelike, then the ruled surface $N$ is said to be of type $N_+^1$.

**iii)** If the central normal vector $\vec{h}$ is timelike and the ruling $\vec{q}$ is spacelike, then the ruled surface $N$ is said to be of type $N_+^2$ [9,19].

The ruled surfaces of type $N_-^1$ and $N_+^1$ are clearly timelike and the ruled surface of type $N_+^2$ is spacelike.

By using these classifications, the parametrization of the ruled surface $N$ can be given as follows,

$$\varphi(s, v) = \vec{k}(s) + v\vec{q}(s), \qquad (4)$$

where $\langle \vec{h}, \vec{h} \rangle = \varepsilon_1 (= \pm 1)$, $\langle \vec{q}, \vec{q} \rangle = \varepsilon_2 (= \pm 1)$.

The set of all bound vectors $\vec{q}(s)$ at the point O constitutes the *directing cone* of the ruled surface $N$. If $\varepsilon_2 = -1$ (resp. $\varepsilon_2 = 1$), the end points of the vectors $\vec{q}(s)$ drive a spherical spacelike (resp. spacelike or timelike) curve $k_1$ on hyperbolic unit sphere $H_0^2$ (resp. on Lorentzian unit sphere $S_1^2$), called the *hyperbolic (resp. Lorentzian) spherical image* of the ruled surface $N$, whose arc is denoted by $s_1$.

For the Frenet vectors $\vec{q}, \vec{h}$ and $\vec{a}$ we have the following Frenet frames of ruled surface $N$:

**i)** If the ruled surface $N$ is timelike ruled surfaces of the type $N_+^1$ or $N_-^1$ then we have

$$\begin{bmatrix} d\vec{q}/ds_1 \\ d\vec{h}/ds_1 \\ d\vec{a}/ds_1 \end{bmatrix} = \begin{bmatrix} 0 & 1 & 0 \\ -\varepsilon_2 & 0 & \kappa \\ 0 & \varepsilon_2\kappa & 0 \end{bmatrix} \begin{bmatrix} \vec{q} \\ \vec{h} \\ \vec{a} \end{bmatrix}. \qquad (5)$$

From (5), the Darboux vector of the Frenet frame $\{O; \vec{q}, \vec{h}, \vec{a}\}$ can be given by $\vec{w}_1 = \varepsilon_2 \kappa \vec{q} - \vec{a}$. Thus, for the derivatives in (5) we can write

$$d\vec{q}/ds_1 = \vec{w}_1 \times \vec{q}, \quad d\vec{h}/ds_1 = \vec{w}_1 \times \vec{h}, \quad d\vec{a}/ds_1 = \vec{w}_1 \times \vec{a},$$

and also we have

$$\vec{q} \times \vec{h} = \varepsilon_2 \vec{a}, \quad \vec{h} \times \vec{a} = -\varepsilon_2 \vec{q}, \quad \vec{a} \times \vec{q} = -\vec{h}. \tag{6}$$

[See 15].

**ii)** If the ruled surface $N$ is spacelike ruled surface then we have

$$\begin{bmatrix} d\vec{q}/ds_1 \\ d\vec{h}/ds_1 \\ d\vec{a}/ds_1 \end{bmatrix} = \begin{bmatrix} 0 & 1 & 0 \\ 1 & 0 & \kappa \\ 0 & \kappa & 0 \end{bmatrix} \begin{bmatrix} \vec{q} \\ \vec{h} \\ \vec{a} \end{bmatrix}. \tag{7}$$

Darboux vector of this frame is $\vec{w}_1 = -\kappa \vec{q} + \vec{a}$. Then the derivatives of the vectors of Frenet frame in (7) can be given by

$$d\vec{q}/ds_1 = \vec{w}_1 \times \vec{q}, \quad d\vec{h}/ds_1 = \vec{w}_1 \times \vec{h}, \quad d\vec{a}/ds_1 = \vec{w}_1 \times \vec{a}$$

and also we have

$$\vec{q} \times \vec{h} = -\vec{a}, \quad \vec{h} \times \vec{a} = -\vec{q}, \quad \vec{a} \times \vec{q} = \vec{h}. \tag{8}$$

[See 26].

In these equations, $s_1$ is the arc of generating curve $k_1$ and $\kappa = \dfrac{ds_3}{ds_1} = \left\| \dfrac{d\vec{a}}{ds} \right\|$ is conical curvature of the directing cone, where $s_3$ is the arc of the spherical curve $k_3$ circumscribed by the bound vector $a$ at the point $O$ [15,16].

## 4. Mannheim Offsets of Timelike Ruled Surfaces in Minkowski 3-space.

Assume that $\varphi$ and $\varphi^*$ be two ruled surfaces in the Minkowski 3-space $IR_1^3$ with the parametrizations

$$\varphi(s,v) = \vec{c}(s) + v\vec{q}(s), \quad \|\vec{q}(s)\| = 1,$$
$$\varphi^*(s,v) = \vec{c}^*(s) + v\vec{q}^*(s), \quad \|\vec{q}^*(s)\| = 1, \tag{9}$$

respectively, where $(\vec{c})$ (resp. $(\vec{c}^*)$) is the striction curve of the ruled surfaces $\varphi$ (resp. $\varphi^*$). Let the Frenet frames of the ruled surfaces $\varphi$ and $\varphi^*$ be $\{\vec{q}, \vec{h}, \vec{a}\}$ and $\{\vec{q}^*, \vec{h}^*, \vec{a}^*\}$, respectively, and $\varphi$ be a timelike ruled surface. The ruled surface $\varphi^*$ is said to be Mannheim offset of the timelike ruled surface $\varphi$ if there exists a one to one correspondence between their rulings such that the asymptotic normal of $\varphi$ is the central normal of $\varphi^*$. In this case, $(\varphi, \varphi^*)$ is called a pair of Mannheim ruled surface. By this definition we can write,

$$\vec{h}^* = \vec{a} \tag{10}$$

and so that by Definition 2.1 and considering the classifications of the ruled surfaces, we have the followings:

**Case 1.** If the ruled surface $\varphi$ is timelike ruled surface of the type $N_-^1$ then by considering (10), Mannheim offset $\varphi^*$ of $\varphi$ is a timelike ruled surface of the type $N_-^1$ or $N_+^1$. If $\varphi$ is of the type $N_-^1$ and $\varphi^*$ is of the type $N_+^1$, we have

$$\begin{pmatrix} \vec{q}^* \\ \vec{h}^* \\ \vec{a}^* \end{pmatrix} = \begin{pmatrix} \sinh\theta & \cosh\theta & 0 \\ 0 & 0 & 1 \\ \cosh\theta & \sinh\theta & 0 \end{pmatrix} \begin{pmatrix} \vec{q} \\ \vec{h} \\ \vec{a} \end{pmatrix}. \qquad (11)$$

Similarly, if $\varphi$ is of the type $N_-^1$ and $\varphi^*$ is of the type $N_-^1$, we have

$$\begin{pmatrix} \vec{q}^* \\ \vec{h}^* \\ \vec{a}^* \end{pmatrix} = \begin{pmatrix} \cosh\theta & \sinh\theta & 0 \\ 0 & 0 & 1 \\ \sinh\theta & \cosh\theta & 0 \end{pmatrix} \begin{pmatrix} \vec{q} \\ \vec{h} \\ \vec{a} \end{pmatrix}. \qquad (12)$$

**Case 2.** If the ruled surface $\varphi$ is timelike ruled surface of the type $N_+^1$, then Mannheim offset $\varphi^*$ of $\varphi$ is a spacelike ruled surface and we have

$$\begin{pmatrix} \vec{q}^* \\ \vec{h}^* \\ \vec{a}^* \end{pmatrix} = \begin{pmatrix} \cos\theta & \sin\theta & 0 \\ 0 & 0 & 1 \\ \sin\theta & -\cos\theta & 0 \end{pmatrix} \begin{pmatrix} \vec{q} \\ \vec{h} \\ \vec{a} \end{pmatrix}. \qquad (13)$$

In (11), (12) and (13), $\theta$ is the angle between the rulings $\vec{q}$ and $\vec{q}^*$.

By definition, the parametrization of $\varphi^*$ can be given by

$$\varphi^*(s,v) = \left[\vec{c}(s) + R(s)\vec{a}(s)\right] + v\vec{q}^*(s). \qquad (14)$$

From the definition of $\vec{h}^*$, we get $\vec{h}^* = \dfrac{d\vec{q}^*}{ds} / \left\|\dfrac{d\vec{q}^*}{ds}\right\|$. So that we have $\dfrac{d\vec{q}^*}{ds} = \lambda \vec{h}^*$, ($\lambda$ is a scalar). Using this equality and the fact that the base curve of $\varphi^*$ is striction curve we get $\left\langle \dfrac{d}{ds}(\vec{c} + R\vec{a}), \vec{a} \right\rangle = 0$. It follows that $\varepsilon_2 \left\|\dfrac{dq}{ds}\right\| \delta_q + \dfrac{dR}{ds} = 0$. Thus we can give the following theorems.

***Theorem 4.1.*** *Let the ruled surface $\varphi^*$ be Mannheim offset of timelike ruled surface $\varphi$ of the type $N_-^1$ or $N_+^1$. Then $\varphi$ is developable timelike ruled surface if and only if $R$ is a constant.*

Now, we can give the characterizations of the Mannheim offsets of a timelike ruled surface according to the classifications of the surface as follows**.**

### 5. Mannheim Offsets of the Timelike Ruled Surfaces of the Type $N_-^1$

Let the ruled surface $\varphi^*$ be Mannheim offset of developable timelike ruled surface $\varphi$ of the type $N_-^1$. By the definition $\varphi^*$ can be of the type $N_+^1$ or $N_-^1$. Then we can give the followings.

***Theorem 5.1.*** *i) Let the timelike ruled surface $\varphi^*$ be a Mannheim offset of the developable timelike ruled surface $\varphi$ of the type $N_-^1$. Then $\varphi^*$ is developable if and only if the following equality holds*

$$\cosh\theta + R\kappa \frac{ds_1}{ds}\sinh\theta = 0, \text{ if } \varphi^* \text{ is of the type } N_+^1, \tag{15}$$

$$\sinh\theta + R\kappa \frac{ds_1}{ds}\cosh\theta = 0, \text{ if } \varphi^* \text{ is of the type } N_-^1. \tag{16}$$

**Proof.** Let the timelike ruled surface $\varphi^*$ of the type $N_+^1$ be developable. Then we have

$$\frac{d\vec{c}^*}{ds} = \mu \vec{q}^*, \tag{17}$$

where $\mu$ is scalar and $s$ is the arc-length parameter of the striction curve $(c)$ of the timelike ruled surface $\varphi$ of the type $N_-^1$. Then from (11) we obtain

$$\frac{d\vec{c}}{ds} + \frac{dR}{ds}\vec{a} + R\frac{ds_1}{ds}\frac{d\vec{a}}{ds_1} = \mu(\sinh\theta \vec{q} + \cosh\theta \vec{h}). \tag{18}$$

From Theorem 4.1 and by using (5) we get

$$\vec{q} - R\frac{ds_1}{ds}\kappa \vec{h} = \mu\sinh\theta \vec{q} + \mu\cosh\theta \vec{h}. \tag{19}$$

From the last equation it follows that

$$\cosh\theta + R\kappa \frac{ds_1}{ds}\sinh\theta = 0. \tag{20}$$

Conversely, if (20) holds then for the tangent vector of the striction curve $(\vec{c}^*)$ of the timelike ruled surface $\varphi^*$ of the type $N_+^1$ we can write

$$\frac{d\vec{c}^*}{ds} = \frac{d}{ds}(\vec{c} + R\vec{a})$$

$$= \vec{q} - R\frac{ds_1}{ds}\kappa \vec{h}$$

$$= \frac{1}{\sinh\theta}(\sinh\theta \vec{q} + \cosh\theta \vec{h})$$

$$= \frac{1}{\sinh\theta}\vec{q}^*$$

Thus $\varphi^*$ is developable.

If the timelike ruled surface $\varphi^*$ is of the type $N_-^1$ and a Mannheim offset of developable timelike ruled surface $\varphi$ of the type $N_-^1$, by making similar calculations it can be easily shown that $\varphi^*$ is developable if and only if following equality holds

$$\sinh\theta + R\kappa \frac{ds_1}{ds}\cosh\theta = 0. \tag{21}$$

**Theorem 5.2.** *Let $\varphi$ be a developable timelike ruled surface of the type $N_-^1$. The developable timelike ruled surface $\varphi^*$ of the type $N_+^1$ or $N_-^1$ is a Mannheim offset of the ruled surface $\varphi$ if and only if the following relationship holds*

$$\frac{d\kappa}{ds} = -\frac{1}{R}\left(R^2\kappa^2\left(\frac{ds_1}{ds}\right)^2 - 1\right) - \frac{1}{ds_1/ds}\frac{d^2 s_1}{ds^2}\kappa. \tag{22}$$

**Proof.** Let developable timelike ruled surface $\varphi^*$ be a Mannheim offset of timelike ruled surface $\varphi$ of the type $N_-^1$ and assume that $\varphi^*$ is of the type $N_+^1$. From Theorem 5.1 we have

$$R\kappa \frac{ds_1}{ds} = -\coth\theta. \tag{23}$$

Using (11) we have

$$\frac{d\vec{q}^*}{ds} = \cosh\theta\left(\frac{d\theta}{ds} + \frac{ds_1}{ds}\right)\vec{q} + \sinh\theta\left(\frac{d\theta}{ds} + \frac{ds_1}{ds}\right)\vec{h} + \cosh\theta\kappa\frac{ds_1}{ds}\vec{a}. \tag{24}$$

From (24) and definition of $\vec{h}^*$ we have

$$\frac{d\theta}{ds} = -\frac{ds_1}{ds}. \tag{25}$$

Differentiating (23) with respect to $s$ and using (25) we get

$$\frac{d\kappa}{ds} = -\frac{1}{R}\left(R^2\kappa^2\left(\frac{ds_1}{ds}\right)^2 - 1\right) - \frac{1}{ds_1/ds}\frac{d^2s_1}{ds^2}\kappa. \tag{26}$$

Conversely, if (26) holds then for nonzero constant scalar $R$ we can define a timelike ruled surface $\varphi^*$ of the type $N_+^1$ as follows

$$\varphi^*(s,v) = c^*(s) + v q^*(s), \tag{27}$$

where $\vec{c}^*(s) = \vec{c}(s) + R\vec{a}(s)$. Since $\varphi^*$ is developable, we have

$$\frac{d\vec{c}^*}{ds} = \frac{ds^*}{ds}\vec{q}^*, \tag{28}$$

where $s$ and $s^*$ are arc-length parameters of striction curves $(\vec{c})$ and $(\vec{c}^*)$, respectively. From (28) we get

$$\frac{ds^*}{ds}\vec{q}^* = \frac{d}{ds}(\vec{c} + R\vec{a}) = \vec{q} - R\kappa\frac{ds_1}{ds}\vec{h}. \tag{29}$$

By taking the derivative of (29) with respect to $s$, we have

$$\frac{d^2s^*}{ds^2}\vec{q}^* + \frac{ds^*}{ds}\frac{d\vec{q}^*}{ds} = -R\kappa\left(\frac{ds_1}{ds}\right)^2\vec{q} + \left(\frac{ds_1}{ds} - R\kappa\frac{d^2s_1}{ds^2} - R\frac{ds_1}{ds}\frac{d\kappa}{ds}\right)\vec{h} - R\kappa^2\left(\frac{ds_1}{ds}\right)^2\vec{a}. \tag{30}$$

From the hypothesis and the definition of $\vec{h}^*$, we get

$$\frac{d^2s^*}{ds^2}\vec{q}^* + \frac{ds^*}{ds}\lambda\vec{h}^* = -R\kappa\left(\frac{ds_1}{ds}\right)^2\vec{q} + R^2\kappa^2\left(\frac{ds_1}{ds}\right)^3\vec{h} - R\kappa^2\left(\frac{ds_1}{ds}\right)^2\vec{a}, \tag{31}$$

where $\lambda$ is a scalar. By taking the vector product of (29) with (31), we obtain

$$\left(\frac{ds^*}{ds}\right)^2\lambda\vec{a}^* = R^2\kappa^3\left(\frac{ds_1}{ds}\right)^3\vec{q} - R\kappa^2\left(\frac{ds_1}{ds}\right)^2\vec{h}. \tag{32}$$

Taking the vector product of (32) with (29), we have

$$-\left(\frac{ds^*}{ds}\right)^3\lambda\vec{h}^* = \left[R^3\kappa^4\left(\frac{ds_1}{ds}\right)^4 - R\kappa^2\left(\frac{ds_1}{ds}\right)^2\right]\vec{a}. \tag{33}$$

It shows that, developable timelike ruled surface $\varphi^*$ of the type $N_+^1$ is a Mannheim offset of timelike ruled surface $\varphi$ of the type $N_-^1$.

If $\varphi^*$ is of the type $N_-^1$, then making the similar calculations it is easily seen that developable timelike ruled surface $\varphi^*$ is a Mannheim offset of timelike ruled surface $\varphi$ of the type $N_-^1$ if and only if the equation (22) holds.

Let now timelike ruled surface $\varphi^*$ be a Mannheim offset of timelike ruled surface $\varphi$ of the type $N_-^1$. If the trajectory ruled surfaces generated by the vectors $\vec{h}^*$ and $\vec{a}^*$ of $\varphi^*$ are denoted by $\varphi_{h^*}$ and $\varphi_{a^*}$, respectively, then we can write

$$\vec{q}_1^* = \vec{h}^* = \vec{a}, \ \vec{h}_1^* = \mp \vec{h}, \ \vec{a}_1^* = \pm \vec{q}, \tag{34}$$

and if $\varphi^*$ is of the type $N_+^1$,

$$\vec{q}_2^* = \vec{a}^* = (\cosh\theta)\vec{q} + (\sinh\theta)\vec{h}, \ \vec{h}_2^* = \mp\vec{a}, \ \vec{a}_2^* = \mp\left((\sinh\theta)\vec{q} + (\cosh\theta)\vec{h}\right), \tag{35}$$

similarly if $\varphi^*$ is of the type $N_-^1$, then

$$\vec{q}_2^* = \vec{a}^* = (\sinh\theta)\vec{q} + (\cosh\theta)\vec{h}, \ \vec{h}_2^* = \mp\vec{a}, \ \vec{a}_2^* = \pm\left((\cosh\theta)\vec{q} + (\sinh\theta)\vec{h}\right), \tag{36}$$

where $\{\vec{q}_1^*, \vec{h}_1^*, \vec{a}_1^*\}$ and $\{\vec{q}_2^*, \vec{h}_2^*, \vec{a}_2^*\}$ are the Frenet Frames of ruled surfaces $\varphi_{h^*}$ and $\varphi_{a^*}$, respectively. Therefore from (34), (35) and (36) we have the following.

**Corollary 5.3.** (a) $\varphi_{h^*}$ *is a Bertrand offset of* $\varphi$.
(b) $\varphi_{a^*}$ *is a Mannheim offset of* $\varphi$.

Now we can give the followings.

Let the timelike ruled surface $\varphi^*$ be a Mannheim offset of the developable timelike ruled surface $\varphi$ of the type $N_-^1$. From (5), (3), (11) and (12), we obtain

$$\delta_{h^*} = -\frac{1}{(ds_1/ds)\kappa}, \tag{37}$$

$$\delta_{a^*} = -\frac{\sinh\theta + R(ds_1/ds)\kappa\cosh\theta}{(ds_1/ds)\kappa\sinh\theta}, \ \text{if } \varphi^* \text{ is of the type } N_+^1. \tag{38}$$

$$\delta_{a^*} = -\frac{\cosh\theta + R(ds_1/ds)\kappa\sinh\theta}{(ds_1/ds)\kappa\cosh\theta}, \ \text{if } \varphi^* \text{ is of the type } N_-^1. \tag{39}$$

Then, we can give the following corollary.

**Corollary 5.4.** (a) $\varphi_{h^*}$ *is not developable while* $\varphi$ *is developable.*
(b) $\varphi_{a^*}$ *is developable while* $\varphi$ *is developable if and only if the following equalities holds,*
$$\sinh\theta + R(ds_1/ds)\kappa\cosh\theta = 0, \ \text{if } \varphi^* \text{ is of the type } N_+^1. \tag{40}$$
$$\cosh\theta + R(ds_1/ds)\kappa\sinh\theta = 0, \ \text{if } \varphi^* \text{ is of the type } N_-^1. \tag{41}$$

**Example 5.1.** Let consider timelike ruled surface $N$ of the type $N_-^1$ given by the parametrization

$$\varphi(s,v) = \left(2\sinh s + v\sqrt{5}\cosh s, \ \sqrt{3}s + 2v, \ 2\cosh s + v\sqrt{5}\cosh s\right)$$
$$= \left(2\sinh s, \ \sqrt{3}s, \ 2\cosh s\right) + v\left(\sqrt{5}\cosh s, \ 2, \ \sqrt{5}\sinh s\right)$$

which is not developable. For the Frenet vectors of the surface we obtain followings

$$\vec{q}(s) = \left(\sqrt{5}\cosh s,\ 2,\ \sqrt{5}\sinh s\right)$$

$$\vec{h}(s) = (\sinh s, 0, \cosh s)$$

$$\vec{a}(s) = \left(-2\cosh s,\ -\sqrt{5},\ -2\sinh s\right)$$

Then the general parametrization of Mannheim offsets of the type $N_+^1$ is given by

$$\varphi^*(s,v) = \left(2\sinh s - 2R(s)\cosh s,\ \sqrt{3}s - \sqrt{5}R(s),\ 2\cosh s - 2R(s)\sinh s\right)$$

$$+ v\left(\sqrt{5}\sinh\theta\cosh s + \cosh\theta\sinh s,\ 2\sinh\theta,\ \sqrt{5}\sinh\theta\sinh s + \cosh\theta\cosh s\right)$$

Since $\varphi$ is not developable, by taking $R(s) = s$, $\theta = \ln(1+\sqrt{2})$ an offset of $\varphi$ is obtained as follows,

$$\varphi^*(s,v) = \left(2\sinh s - 2s\cosh s,\ \sqrt{3}s - \sqrt{5}s,\ 2\cosh s - 2s\sinh s\right)$$

$$+ v\left(\sqrt{5}\cosh s + \sqrt{2}\sinh s,\ 2,\ \sqrt{5}\sinh s + \sqrt{2}\cosh s\right)$$

Similarly, the general parametrization of Mannheim offsets of the type $N_-^1$ is

$$\varphi^*(s,v) = \left(2\sinh s - 2R(s)\cosh s,\ \sqrt{3}s - \sqrt{5}R(s),\ 2\cosh s - 2R(s)\sinh s\right)$$

$$+ v\left(\sqrt{5}\cosh\theta\cosh s + \sinh\theta\sinh s,\ 2\cosh\theta,\ \sqrt{5}\cosh\theta\sinh s + \sinh\theta\cosh s\right)$$

Then by considering special case $R(s) = 2s$, $\theta = \dfrac{s}{4}$ a Mannheim offset of the type $N_-^1$ is obtained as follows

$$\varphi^*(s,v) = \left(2\sinh s - 4s\cosh s,\ \sqrt{3}s - 2\sqrt{5}s,\ 2\cosh s - 4s\sinh s\right)$$

$$+ v\left(\sqrt{5}\cosh\left(\frac{s}{4}\right)\cosh s + \sinh\left(\frac{s}{4}\right)\sinh s,\ 2\cosh\left(\frac{s}{4}\right),\right.$$

$$\left.\sqrt{5}\cosh\left(\frac{s}{4}\right)\sinh s + \sinh\left(\frac{s}{4}\right)\cosh s\right)$$

## 6. Mannheim Offsets of the Timelike Ruled Surfaces of the Type $N_+^1$

Let the ruled surface $\varphi^*$ be Mannheim offset of the developable timelike ruled surface $\varphi$ of the type $N_+^1$. By the definition, $\varphi^*$ is a spacelike ruled surface. Then we can give the following theorems and corollaries. The proofs of those can be given by the same ways as given in Section 5.

**Theorem 6.1.** *Let the spacelike ruled surface $\varphi^*$ be a Mannheim offset of the developable timelike ruled surface $\varphi$ of the type $N_+^1$. Then $\varphi^*$ is developable if and only if the following equality holds*

$$\sin\theta - R\kappa\frac{ds_1}{ds}\cos\theta = 0. \tag{42}$$

**Theorem 6.2.** *Let $\varphi$ be a developable timelike ruled surface of the type $N_+^1$. The developable spacelike ruled surface $\varphi^*$ is a Mannheim offset of the ruled surface $\varphi$ if and only if the following relationship holds*

$$\frac{d\kappa}{ds} = -\frac{1}{R}\left(R^2\kappa^2\left(\frac{ds_1}{ds}\right)^2 + 1\right) - \frac{1}{ds_1/ds}\frac{d^2s_1}{ds^2}\kappa. \tag{43}$$

Let now the timelike ruled surface $\varphi^*$ be a Mannheim offset of the timelike ruled surface $\varphi$ of the type $N_+^1$. If the trajectory ruled surfaces generated by the vectors $\vec{h}^*$ and $\vec{a}^*$ of $\varphi^*$ are denoted by $\varphi_{h^*}$ and $\varphi_{a^*}$, respectively, then we can write

$$\vec{q}_1^* = \vec{h}^* = \vec{a},\ \vec{h}_1^* = \mp\vec{h},\ \vec{a}_1^* = \mp\vec{q} \tag{44}$$

$$\vec{q}_2^* = \vec{a}^* = (\sin\theta)\vec{q} - (\cos\theta)\vec{h},\ \vec{h}_2^* = \mp\vec{a},\ \vec{a}_2^* = \mp\left((\cos\theta)\vec{q} + (\sin\theta)\vec{h}\right), \tag{45}$$

where $\{\vec{q}_1^*, \vec{h}_1^*, \vec{a}_1^*\}$ and $\{\vec{q}_2^*, \vec{h}_2^*, \vec{a}_2^*\}$ are the Frenet Frames of ruled surfaces $\varphi_{h^*}$ and $\varphi_{a^*}$, respectively. Therefore from (44) and (45) we have the following corollary:

**Corollary 6.3.** (a) $\varphi_{h^*}$ is a Bertrand offset of $\varphi$.

(b) $\varphi_{a^*}$ is a Mannheim offset of $\varphi$.

Now we can give the followings.

Let the spacelike ruled surface $\varphi^*$ be a Mannheim offset of developable timelike ruled surface $\varphi$ of the type $N_+^1$. From (3), (5) and (13), we obtain

$$\delta_{h^*} = \frac{1}{(ds_1/ds)\kappa}, \tag{46}$$

$$\delta_{a^*} = \frac{\cos\theta + R(ds_1/ds)\kappa\sin\theta}{(ds_1/ds)\kappa\cos\theta}. \tag{47}$$

Then, we can give the following corollary.

**Corollary 6.4.** (a) $\varphi_{h^*}$ is not developable while $\varphi$ is developable.

(b) $\varphi_{a^*}$ is developable while $\varphi$ is developable if and only if the following equality holds,

$$\cos\theta + R(ds_1/ds)\kappa\sin\theta = 0. \tag{48}$$

**Example 6.1.** Let consider timelike ruled surface of the type $N_+^1$ given by the parametrization

$$\varphi(s,v) = \left(s,\ \sqrt{2}\sin s,\ -\sqrt{2}\cos s\right) + v\left(0,\ \cos s,\ \sin s\right)$$

which is not developable. For this surface we obtain the following Frenet vectors

$\vec{q}(s) = (0,\ \cos s,\ \sin s)$

$\vec{h}(s) = (0,\ -\sin s, \cos s)$

$\vec{a}(s) = (1,0,0)$.

Then the general parametrization of spacelike Mannheim offsets of the type is

$$\varphi^*(s,v) = \left(s + R(s), \sqrt{2}\sin s,\ -\sqrt{2}\cos s\right)$$

$$+ v\left(0,\ \cos\theta\cos s - \sin\theta\sin s,\ \cos\theta\sin s + \sin\theta\cos s\right)$$

By taking $R(s) = s^2$, $\theta = \pi/4$ an offset of $\varphi$ is obtained as follows,

$$\varphi^*(s,v) = \left(s+s^2,\ \sqrt{2}\sin s,\ -\sqrt{2}\cos s\right)$$
$$+ v\left(0,\ \frac{\sqrt{2}}{2}\cos s - \frac{\sqrt{2}}{2}\sin s,\ \frac{\sqrt{2}}{2}\sin s + \frac{\sqrt{2}}{2}\cos s\right).$$

**7. Conclusions**

In this paper, Mannheim offsets of the timelike ruled surfaces have been developed in Minkowski 3-space $IR_1^3$. It is shown that according to the classifications of the ruled surfaces in Minkowski 3-space $IR_1^3$, the Mannheim offsets of a timelike ruled surface may be timelike or spacelike. Furthermore, developable timelike ruled surfaces can have a developable timelike or spacelike Mannheim offset if the derivative of the conical curvature $\kappa$ of the directing cone holds an equation given by (22) or (43). Since the offset surfaces have an important role in the study of design problems, these obtained results can be let the new studies on such problems.

**References**
[1] Abdel-All, N.H., Abdel-Baky, R.A., Hamdoon, F.M., Ruled surfaces with timelike rulings, App. Math. and Comp., 147 (2004) 241–253.
[2] Beem, J.K., Ehrlich, P.E., Global Lorentzian Geometry, Marcel Dekker, New York, (1981).
[3] Birman, G., Nomizo, K., Trigonometry in Lorentzian Geometry, Ann. Math. Mont., 91, (9) (1984) 543-549.
[4] Dillen, F., Sodsiri, W., Ruled surfaces of Weingarten type in Minkowski 3-space, J. Geom., 83, (2005) 10-21.
[5] Karger, A., Novak, J., Space Kinematics and Lie Groups, STNL Publishers of Technical Lit., Prague, Czechoslovakia (1978).
[6] Kim, Y.H., Yoon, W.D., Classification of ruled surfaces in Minkowski 3-space, J. of Geom. and Physics, 49 (2004) 89-100.
[7] Küçük, A., On the developable timelike trajectory ruled surfaces in Lorentz 3-space $IR_1^3$, App. Math. and Comp., 157 (2004) 483-489.
[8] Liu, H., Wang, F., Mannheim partner curves in 3-space, Journal of Geometry, vol. 88, no. 1-2 (2008) 120-126.
[9] O'Neill, B., Semi-Riemannian Geometry with Applications to Relativity, Academic Press, London, (1983).
[10] Orbay, K., Kasap, E., Aydemir, İ, Mannheim Offsets of Ruled Surfaces, Mathematical Problems in Engineering, Volume 2009, Article ID 160917.
[11] Peternel, M., Pottmann, H., Ravani, B., On the computational geometry of ruled surfaces, Comp. Aided Geom. Design, 31 (1999) 17-32.
[12] Pottmann, H., Lü, W., Ravani, B., Rational ruled surfaces and their offsets, Graphical Models and Image Procesing, vol. 58, no. 6 (1996) 544-552.
[13] Ravani, B., Ku, T.S., Bertrand Offsets of ruled and developable surfaces, Comp. Aided Geom. Design, (23), No. 2, (1991).
[14] Turgut, A., Hacısalihoğlu H.H, Timelike ruled surfaces in the Minkowski 3-space, Far East J. Math. Sci. 5 (1) (1997) 83–90
[15] Uğurlu H.H., Önder M., Instantaneous Rotation vectors of Skew Timelike Ruled Surfaces in Minkowski 3-space, VI. Geometry Symposium, 01-04 July 2008, Bursa, Turkey.

[16] Uğurlu, H.H., Önder, M., On Frenet Frames and Frenet Invariants of Skew Spacelike Ruled Surfaces in Minkowski 3-space, VII. Geometry Symposium, 7-10 July 2009, Kırşehir, Turkey.

[17] Uğurlu, H.H., Çalışkan, A., Darboux Ani Dönme Vektörleri ile Spacelike ve Timelike Yüzeyler Geometrisi, Celal Bayar Üniversitesi Yayınları, Yayın no: 0006, 2012.

[18] Wang, F., Liu, H., Mannheim partner curves in 3-Euclidean space, Mathematics in Practice and Theory, vol. 37, no. 1, pp. 141-143, 2007.